\documentclass[letterpaper, 10 pt, conference]{ieeeconf}

\IEEEoverridecommandlockouts                              
\title{\LARGE \bf SuperSCS: fast and accurate large-scale conic optimization}

\author{%
    Pantelis~Sopasakis$^{*}$,
    Krina~Menounou$^{\dagger}$
    and Panagiotis~Patrinos$^{\dagger}$%
\thanks{$^{*}$ Queen's University Belfast, School of Electronics, Electrical Engineering and Computer Science,
               Centre for Intelligent Autonomous Manufacturing Systems, BT9 5AH, Northern Ireland, UK.
               Email: \texttt{p.sopasakis@qub.ac.uk}}
\thanks{$^{\dagger}$ KU Leuven, Department of Electrical Engineering (ESAT), STADIUS, Kasteelpark 10, 3001 Leuven, Belgium.}    
\thanks{This work is accompanied by the open-source (licensed with the MIT licence) 
        free software \texttt{SuperSCS} which 
        is available online at \href{https://kul-forbes.github.io/scs/}{https://kul-forbes.github.io/scs/}.
        }
\thanks{P. Sopasakis was supported by European Union's Horizon 2020 research and 
        innovation programme (KIOS CoE) under Grant No. 739551.
        The work of the third author was supported by: FWO projects: No. G086318N; No. G086518N; 
        Fonds de la Recherche Scientifique --FNRS, the Fonds Wetenschappelijk Onderzoek 
        -- Vlaanderen under EOS Project No. 30468160 (SeLMA) and 
        Research Council KU Leuven C1 project No. C14/18/068.}
    }

\usepackage{nicefrac}
\usepackage{units}
\usepackage{mathematix}
\usepackage{acronym,etoolbox}
\usepackage{wrapfig}
\usepackage{hyperref}
\usepackage[font=footnotesize,labelfont=sc]{caption}

\input{acronym-config.tex}

\acrodef{dm}[DM]{Dolan-Mor{\'e}}
\acrodef{aa}[AA]{Anderson's acceleration}
\acrodef{rb}[RB]{restarted Broyden}
\acrodef{hsde}[HSDE]{homogeneous self-dual embedding}
\acrodef{scs}[SCS]{splitting cone solver}
\acrodef{drs}[DRS]{Douglas-Rachford splitting}
\acrodef{sdp}[SDP]{semidefinite program}

\hypersetup{
    unicode=false,          
    pdftoolbar=true,        
    pdfmenubar=true,        
    pdffitwindow=false,     
    pdfstartview={FitH},    
    pdfnewwindow=true,      
    colorlinks=true,        
    linkcolor=red,   
    citecolor=blue,         
    filecolor=magenta,      
    urlcolor=cyan           
}

\newcommand{\K}{\mathcal{K}}
\newcommand{\tttop}{{\scriptscriptstyle \top}}
\newcommand{\ttplus}{{\scriptscriptstyle +}}
\newcommand{\ttperp}{{\scriptscriptstyle \perp}}

\begin{document}

\maketitle

\begin{abstract}
We present SuperSCS: a fast and accurate method for solving large-scale convex conic problems. SuperSCS combines the SuperMann algorithmic framework with the Douglas-Rachford splitting which is applied on the homogeneous self-dual embedding of conic optimization problems: a model for conic optimization problems which simultaneously encodes the optimality conditions and infeasibility/unboundedness certificates for the original problem. SuperMann allows the use of fast quasi-Newtonian directions such as a modified restarted Broyden-type direction and Anderson's acceleration.
\end{abstract}

\section{Introduction}

Conic optimization problems are of central importance in convex optimization
as several solvers and parsers such as CVX \cite{cvx}, CVXPy \cite{cvxpy}, 
YALMIP \cite{yalmip} and MOSEK \cite{mosek} transform given 
problems into a conic representation. Indeed, all convex optimization problems
can be cast in the standard form of a conic optimization problem.

Various interior point methods have been proposed for solving conic optimization
problems~\cite{S98guide,sdpt3,ecos2013}.
Interior point methods can achieve high accuracy, yet, do not scale well with 
the problem size. 
On the other hand, first-order methods have low per-iteration 
cost and minimum memory requirements, therefore, are better suited for large-scale
problems~\cite{ODonoghue2016}. 
Recent research has turned to first-order methods for large-scale 
conic problems such as SDPs \cite{ZhengFantuzzi2017}.
However, their convergence rate is at most Q-linear
with a Q-factor often close to one, especially for ill-conditioned problems.

The KKT conditions of a conic optimization problem together with conditions 
for the detection of infeasibility or unboundedness can be combined in 
a convex feasibility problem known as the \ac{hsde} \cite{hsdeOriginal}. 
The \ac{hsde} has been used in both interior point \cite{S98guide} 
and first-order methods \cite{ODonoghue2016}.

In this paper we present a numerical optimization method for solving the 
\ac{hsde}. The \ac{hsde} is first cast as a variational inequality which 
can be equivalently seen as a monotone inclusion. We observe that the 
\ac{scs} presented in \cite{ODonoghue2016} can be interpreted as the application
of the \ac{drs} to that monotone inclusion. We then apply the reverse 
splitting to that monotone inclusion --- which is a firmly nonexpansive 
operator --- and employ the SuperMann scheme \cite{supermann} which allows the 
use of quasi-Newtonian directions such as restarted Broyden directions and 
Anderson's acceleration \cite{WalkerNiAnderson2011,Fang2009}. We call the 
resulting method \textit{SuperSCS}.
This way, SuperSCS can achieve a fast convergence rate while retaining a low 
per-iteration cost. In fact, SuperSCS uses exactly the same oracle as \ac{scs}.



\section{Mathematical preliminaries}
We denote by \(\R\), \(\R_{\ttplus}\), \(\R^n\) and \(\R^{m\times n}\) the sets of 
real numbers, non-negative reals, \(n\)-dimensional 
real vectors and \(m\)-by-\(n\) real matrices respectively. We denote the transpose of a 
matrix \(A\) by \(A^{\tttop}\). For two vectors \(x, y\in \R^n\), we denote by \(\<x, y\> = x^{\tttop}y\)
their standard inner product.
Let \(E\) be a vector space in \(\R^n\). We define the \textit{orthogonal complement}
of \(E\) in \(\R^n\) to be the vector space \(E^{\ttperp} = \{y\in\R^n {}\mid{} \<y, x\> = 0, \forall x\in E\}\).

A set \(\K\subseteq\R^n\) is called a \textit{convex cone} if it is convex and
\(\lambda x {}\in{} \K\) for every \(x{}\in{}\K\) and $\lambda>0$.
The binary relation \(x \succcurlyeq_{\K}  y\) is interpreted as \(x-y \in \K\).
The \textit{dual cone} of \(\K\) is defined as 
\(
 \K^* = \{x^* {}\mid{} \<x^*, x\> \geq 0, \forall x\in\K\}.
\)

A few examples of convex cones of interest are:
(i) the zero cone \(\K_n^{\mathrm{f}}=\{0\}^n\),
(ii) the cone of symmetric positive semidefinite matrices 
       \(
			\K^{\mathrm{s}}_n
	{}={}
			\{
			  x\in\R^{n(n+1)/2} 			  
			  {}\mid{}
			  \operatorname{mat}(x): \text{pos. definite}
			\},
       \)
       where \(\operatorname{mat}:\R^{n(n+1)/2} \to \R^{n\times n}\) is defined by
       \[
        \operatorname{mat}(x) = \tfrac{1}{\sqrt{2}}\smallmat{
        \sqrt{2}x_1 & x_2 & \cdots & x_n
        \\
        x_2 & \sqrt{2}x_{n+1} & \cdots & x_{2n-1}
        \\
        \vdots & \vdots & \ddots & \vdots
        \\
        x_n & x_{2n-1} & \cdots & \sqrt{2}x_{n(n+1)/2}
        },
       \]
(iii) the second-order cone 
	\(
			\K^{\mathrm{q}}_n
	  {}={}
			\{
			      z=(x,t)
			      {}:{} 
			      x\in\R^{n-1}, 
			      t{}\in{}\R
			  {}\mid{} 
			      \|x\|_2\leq t
			\},
	\)
(iv) the positive orthant 
	\(
			\K^{\mathrm{l}}_n
	{}={}
			\{
			      x\in\R^n
			  {}\mid{}
			      x\geq 0
			\}
	\) and
(v) the three-dimensional exponential cone, 
       \(
	      \K^{\mathrm{e}}
	{}={}
	      \operatorname{cl}
	      \{
		(x_1,x_2,x_3) 
		{}\mid{}
		x_1{}\geq{}x_2 e^{x_3/x_2},{}
		x_2>0\}
	.\) 

The \textit{normal cone} of a nonempty closed convex set \(C\)
is the set-valued mapping \(N_C(x) = \{g {}\mid{} \<g, y-x\> \leq 0, \forall 
y \in C\}\) for \(x\in C\) and \(N_C(x) = \emptyset\) for \(x\notin C\). 
The \emph{Euclidean projection} of $x$ on $C$ is denoted by $\Pi_C$.

\section{Conic programs}\label{sec:conic_programs}
A cone program is an optimization problem of the form
\begin{equation}\tag{$\mathcal{P}$}\label{eq:cone-program-1}
\begin{array}{rl}
\displaystyle{\minimize_{x\in\R^n}}&\ \<c, x\>\\
\stt&\ b{}-{}Ax {}={} s,\ {} s{}\in{}\K,
\end{array}
\end{equation}
where \(A\in\R^{m\times n}\) is a possibly sparse matrix and \(\K\) is a nonempty 
closed convex cone.
  
The vast majority of convex optimization problems of practical interest can be 
represented in the above form~\cite{ODonoghue2016,BenTal2001}.  
Indeed, cone programs can be thought of as 
a universal representation for all convex problems of practical interest and 
many convex optimization solvers first transform the given problem into this form.

The dual of~\eqref{eq:cone-program-1} is given by 
\cite[Sec. 1.4.3]{BenTal2001}
\begin{equation}\tag{$\mathcal{D}$}\label{eq:conic-dual}
\begin{array}{rl}
\displaystyle{\minimize_{y\in\R^m}}&\ \<b,y\>\\
\stt{}		&\ y\in\K^*,\
						 {} A^{\tttop}y {}+{} c {}={} 0
\end{array}
\end{equation}
Let \(p^\star\) be the optimal value of~\eqref{eq:cone-program-1} and \(d^\star\)
be the optimal value of~\eqref{eq:conic-dual}. Strong duality holds 
(\(p^\star = -d^{\star}\)) if the primal or the dual problem are strictly 
feasible \cite[Thm. 1.4.2]{BenTal2001}.

Whenever strong duality holds, the following KKT conditions are necessary and 
sufficient for optimality of  
\((x^{\star}, s^{\star}, y^{\star})\in\R^{n}\times \R^{m} \times \R^{m}\): 
\begin{gather}
 \begin{aligned}
   Ax^\star + s^\star {}={} b,
   \ 
   s^\star \in \K,
   \
   y^\star \in \K^*,
   \\
   A^{\tttop} y^\star + c = 0,
   \
   \<y^\star,s^\star\> = 0.
 \end{aligned}
\end{gather} 
%

Infeasibility and unboundedness conditions are provided 
by the so-called \textit{theorems of the alternative}. The \textit{weak} 
theorems of  the alternative state that 1) Either primal feasibility holds, or
there is a \(y\) with 
\(A^{\tttop} y = 0\), \(y \succeq_{\K^*} 0\) and \(\<b,y\><0\), 
and, similarly, 2) Either dual feasibility holds or 
there is a \(x\) so that \(Ax\succeq_{\K} 0\) and  \(\<c,x\>\leq 0\) 
\cite[Sec. 1.4.7]{BenTal2001}.

\subsection{Homogeneous self-dual embedding}
In this section we present a key result which is due to Ye \textit{et al.}%
~\cite{hsdeOriginal}: the \ac{hsde}, which is a 
feasibility problem which simultaneously describes the optimality, 
(in)feasibility and (un)boundedness conditions of a conic optimization problem. 
Solving the \ac{hsde} yields a solution of the original conic optimization problem, 
when one exists, or a certificate of infeasibility or unboundedness.
%
%
We start by considering the following feasibility problem in \((\chi, \varsigma, \psi, \tau, \kappa)\)
\begin{subequations}
\begin{align}\label{eq:hsde}
 \smallmat{
  0
  \\
  \varsigma
  \\
  \kappa
 }
{}={}
 Q
 \smallmat{
  \chi
  \\
  \psi
  \\
  \tau
 },
 \ 
 \varsigma\in\K,\ 
 \psi\in\K^*,\ 
 \tau\geq 0,\ 
 \kappa \geq 0,
\end{align}
where
\begin{equation}
	Q
{}\dfn{}
	\smallmat{
	   0   &  A^* & c     \\
	  -A   &  0   & b     \\
	  -c^* & -b^* & 0
	}
\end{equation}
\end{subequations}
Note that for \(\tau^{\star}=1\) and \(\kappa^{\star}=0\), the above equations reduce to the 
primal-dual optimality conditions. As shown in \cite{hsdeOriginal}, the solutions of 
\eqref{eq:hsde} satisfy 
\(
	\kappa^{\star}\tau^{\star}
{}={}
	0
\),
\emph{i.e.}, at least one of \(\kappa^\star\) and \(\tau^\star\) must be zero.
In particular, if \(\kappa=0\) and \(\tau>0\), then the triplet \((x^{\star}, y^{\star}, s^{\star})\) with
\begin{align*}
	x^{\star}
{}={}
	\chi^{\star}/\tau^{\star},
\
	y^{\star}
{}={}
	\psi^{\star}/\tau^{\star},
\
	s^\star
{}={}
	\varsigma^{\star}/\tau^{\star},
\end{align*}
is a primal-dual solution of~\eqref{eq:cone-program-1} and~\eqref{eq:conic-dual}.
If instead \(\tau^{\star}=0\) and \(\kappa>0\), then the problem is either primal- or dual-infeasible.
If \(\tau{}={}\kappa{}={}0\), no conclusion can be drawn.

We define \(u{}={}(\chi,\psi,\tau)\) and \(v=(0,\varsigma,\kappa)\). 
The self-dual embedding reduces to the problem of determining 
\(u\) and \(v\) such that \(Qu=v\) with 
\(
(u,v)\in\mathcal{C}\times\mathcal{C}^*
\) 
where
\(
\mathcal{C} {}\dfn{} \R^n{}\times{}\K^*{}\times{}\R_+.
\)
This is equivalent to the variational inequality
\begin{align}%
 \label{eq:monotone_inclusion}
 0\in Qu + N_{\mathcal{C}}(u),
\end{align}%
Indeed, for all \(u\in\mathcal{C}\), 
\begin{align*}
	N_{\mathcal{C}}(u) 
{}={}& 
	\{ 
	      y 
	{}\mid{} 
	      \<v-u, y\> \leq 0, \forall v \in \mathcal{C}
	\}
\\
{}={}&
	\{u\}^{\ttperp}
{}\cap{}
	\{ 
	      y 
	{}\mid{} 
	      \<v, y\> \leq 0, \forall v \in \mathcal{C}
	\}
{}={}
	\{u\}^{\ttperp} \cap (-\mathcal{C}^*),
\end{align*}
where the second equality follows by considering \(v = \nicefrac{1}{2}u\) and 
\(v=\nicefrac{3}{2}u\), which both belong to the cone \(\mathcal{C}\) 
(see also \cite[Ex.~5.2.6(a)]{hiriartBook}). From that and the fact that 
\(Qu \in \{u\}^{\ttperp}\) since \(Q\) is skew-symmetric, the equivalence of 
the \ac{hsde} in \eqref{eq:hsde} and the variational inequality in 
\eqref{eq:monotone_inclusion} follows.
%
Equation \eqref{eq:monotone_inclusion} is a monotone inclusion which 
can be solved using operator theory machinery as we discuss in the following section.
\section{SuperSCS}\label{sec:superscs}

\subsection{SCS and DRS}
Since \(Q\) is a skew-symmetric linear operator, it is maximally monotone.
Being the normal cone of a convex set, \(N_{\mathcal{C}}\) is maximally monotone as well. Additionally, because of~\cite[Cor. 24.4(i)]{bauschkeBook},
\(Q+N_{\mathcal{C}}\) is maximally monotone. Therefore, we may apply the Douglas-Rachford splitting
on the monotone inclusion~\eqref{eq:monotone_inclusion}. 
The SCS algorithm~\cite{ODonoghue2016} is precisely the application of 
the \ac{drs} to \(N_{\mathcal{C}}{}+{}Q\) leading to the 
iterations discussed in~\cite[Sec. 7.3]{monotonePrimer}. This observation furnishes a short and elegant interpretation of SCS.
Here, on the other hand, we consider the reverse splitting, \(Q+N_{\mathcal{C}}\), 
which leads to the following \ac{drs} iterations
\begin{subequations}
  \label{eq:drs}
  \begin{align}
	    \tilde{u}^{\nu} 
  {}={}& 
	   (I+Q)^{-1}(u^{\nu})
  \label{eq:drs-linear-system}
  \\
	    \bar{u}^{\nu} 
  {}={}& 
	     \Pi_{\mathcal{C}}(2\tilde{u}^{\nu} - u^{\nu})
  \label{eq:drs-projection}
  \\
	    u^{\nu+1} 
  {}={}& 
	    u^{\nu} {}+{} \bar{u}^{\nu} {}-{} \tilde{u}^{\nu}.
  \label{eq:drs-last}	    
 \end{align}
\end{subequations}
For any initial guess \(u^0\), the iterates \(u^\nu\) converge to a point \(u^\star\)
which satisfies the monotone inclusion~\eqref{eq:monotone_inclusion}~\cite[Thm. 25.6(i), (iv)]{bauschkeBook}.
The linear system in~\eqref{eq:drs-linear-system} can be either solved ``directly'' 
using a sparse LDL factorization or ``indirectly'' by means of the conjugate gradient 
method~\cite{ODonoghue2016}. The projection on \(\mathcal{C}\) in~\eqref{eq:drs-projection}
essentially requires that we be able to project on \(\mathcal{K}^*\).

The iterative method~\eqref{eq:drs} can be concisely written as 
\begin{align}\label{eq:T-updates}
 u^{\nu+1} = Tu^\nu,
\end{align}
where $T:\R^{N}\to\R^{N}$ is given by
\(
Tu=u+\Pi_{\mathcal{C}}(2(I+Q)^{-1}u-u)-(I+Q)^{-1}u
\)
and is firmly nonexpansive~\cite[Chap. 26]{bauschkeBook}.
As such it fits the Krasnosel'skii-Mann framework~\cite[Sec.~5.2]{bauschkeBook} leading
to the relaxed iterations
\begin{equation}\label{eq:km-updates}
 u^{\nu+1} {}={} (1-\lambda)u^\nu + \lambda Tu^{\nu},
\end{equation}
with \(\lambda \in (0,2)\) and, as a result,
it fits the SuperMann framework~\cite{supermann}.

\subsection{SuperSCS: SuperMann meets SCS}
SuperMann considers the problem of finding a fixed-point 
\(x^\star{}\in{}\fix T\) from the viewpoint of finding a zero of the residual 
operator
\begin{equation}
 R = I-T.
\end{equation}
SuperMann, instead of applying Krasnosel'skii-Mann-type updates of the 
form~\eqref{eq:km-updates}, takes extragradient-type updates of the general form
\begin{subequations}
 \begin{align}
	  w^{\nu} 
  {}={}&
	  u^{\nu} {}+{} \alpha_\nu d^\nu,
  \\
	 u^{\nu+1}
  {}={}&
	 u^\nu - \zeta_{\nu} Rw^\nu,\label{eq:K2step}
 \end{align}
\end{subequations}
where 
\(
	d^{\nu}
\)
are fast, e.g., quasi-Newtonian, directions and scalar parameters \(\alpha_{\nu}\)
and \(\zeta_{\nu}\) are appropriately chosen so as to guarantee global convergence.

At each step we perform backtracking line search on $\alpha_\nu$ until we either 
trigger fast convergence (K1 steps) or ensure global convergence (K2 steps) as 
shown in Algorithm~\ref{alg:supermann}. The K2 step, cf.~\eqref{eq:K2step}, 
can be interpreted as a projection of the current iterate on a hyperplane 
generated by $w^\nu$, separating the set of fixed points from $u^\nu$, and thus 
guarantees that every iterate comes closer to fixed point set.
Alongside, a sufficient decrease of the norm of the residual, \(\|Ru^{\nu}\|\),
may trigger a ``blind update'' (K0 steps) of the form 
\(
	u^{\nu+1} 
{}={} 
	u^{\nu} 
{}+{} 
	d^{\nu},
\)
where no line search iterations need to be executed.
 \renewcommand{\algorithmicrequire}{\textbf{Input:}}
 \renewcommand{\algorithmicensure}{\textbf{Output:}}
\begin{algorithm} 
\footnotesize
 \begin{algorithmic}
 \REQUIRE \(c_0, c_1, q\in [0,1)\), \(\sigma\in (0,1)\), \(u^0\), \(\lambda\in(0,2)\) and \(\epsilon>0\)%
 \STATE \(\eta_0{}\gets{}\|Ru^0\|\), \(r_{\text{safe}}\gets \eta_0\)
 \FOR{\(\nu=0,1,\ldots\)}
      \STATE Check termination with tolerance \(\epsilon\) (Sec.~\ref{sec:termination})
      \STATE Choose direction \(d^{\nu}\) (Sec.~\ref{sec:quasi-newtonian}), let \(\alpha_{\nu}{}\gets{}1\)
      \IF {
		\(
			\|Ru^\nu\| 
		  {}\leq{} 
			c_0 \eta_\nu
		\)
	      }
	      \STATE 
		  (K0)
		  \(
			u^{\nu+1}
		    {}\gets{}
			w^\nu
		  \),
		  \(
		      \eta^{\nu+1}
		    {}\gets{}
		      \|Ru^{\nu}\|
		  \)
\ELSE
      \STATE
	    \(
		    \eta_{\nu+1}
	      {}\gets{}
		    \eta_{\nu}
	    \)
      \LOOP
	  \STATE 
	      \(
			w^\nu
		{}\gets{} 
			u^\nu 
		{}+{}
			\alpha_\nu d^\nu
	      \)
	  and
	      \(
			\rho_{\nu}
		{}\gets{}
			\langle Rw^\nu, u^\nu {-} Tw^\nu\rangle 
	      \)
	  \IF {
		\(
		      \|Ru^\nu\|
		  {}\leq{}
		      r_{\text{safe}}
		\)
		and
		\(
		      \|Rw^\nu\|
		  {}\leq{} 
		      c_1 \|Ru^\nu\|
		\)
	      }
	    \STATE 
		  (K1) 
		  \(
			    u^{\nu+1} 
		    \,{\gets} \,
			    w^\nu
		  \),
		  \(
			  r_{\text{safe}}
		    {}\gets{}
			  \|Rw^\nu\| + q^\nu \eta_0
		  \),
		   exit loop
	  \ELSIF {
		  \(
		      \rho_{\nu}		  
		   {}\geq{} 
		      \sigma \|Ru^\nu\| \|Rw^\nu\|
		  \)
		 } 
	    \STATE (K2) 
		  \(
			  u^{\nu+1} 
		   {}\gets{} 
			  u^\nu 
		   {}-{} 
			  \lambda\tfrac{\rho_\nu}{\|Rw^\nu\|^2}Rw^\nu
		  \) and exit loop
	  \ELSE
	      \STATE \(\alpha_{\nu}{}\gets{}\nicefrac{\alpha_{\nu}}{2}\)
	  \ENDIF	  
      \ENDLOOP
   \ENDIF
 \ENDFOR{} 
\end{algorithmic} 
\caption{SuperSCS algorithm}\label{alg:supermann}
\end{algorithm}

By exploiting the structure of \(T\) and, in particular, linearity of \((I+Q)^{-1}\), 
we may avoid evaluating linear system solves at every backtracking
step. Instead, we only need to evaluate \(\Pi_{\mathcal{C}}\)
once in every backtracking iteration. In particular, for \(w^\nu = u^\nu + \alpha d^\nu\), we have 
\begin{align*}
      \tilde{w}^\nu
{}={}&
      (I+Q)^{-1}w^\nu
\\
{}={}&
      (I+Q)^{-1}u^\nu {}+{} \alpha (I+Q)^{-1}d^\nu 
{}={} 
      \tilde{u}^\nu+\alpha\tilde{d}^\nu
\end{align*}
where $\tilde{u}^\nu$ has already been computed, since it is needed in the 
evaluation of $Ru^\nu$, while \(\tilde{d}^\nu\) solves \((I+Q)\tilde{d}^\nu = d^\nu\). 
The computation of \(\tilde{d}^\nu\), which is the most costly operation, is 
performed only once, before the backtracking procedure takes place.
The computation of the fixed-point residual of \(w\) is also easily computed 
by \(Rw^{\nu} {}={} \tilde{w}^{\nu}-\proj_{\mathcal{C}}(2\tilde{w}^{\nu}-w^{\nu})\).

%

Overall, save the computation of the residuals, at every iteration of SuperSCS 
we need to solve the linear system \eqref{eq:drs-linear-system} twice and 
invoke \(\Pi_{\mathcal{C}}\) exactly \(1+l_\nu\) times, where \(l_{\nu}\) is the number of
backtracks.

\subsection{Termination}\label{sec:termination}
\begin{subequations}
The algorithm is terminated when an approximate optimal solution is found based on its
relative primal and dual residuals and relative duality gap, provided such a solution exists.
At iteration \(\nu\) let 
\(
	{u}^\nu 
  {}={}
	({\chi}^{\nu}, {\psi}^{\nu}, {\tau}^{\nu})
\)%
,
\(
	\bar{u}^\nu 
  {}={}
	(\bar{\chi}^{\nu}, \bar{\psi}^{\nu}, \bar{\tau}^{\nu})
\)
and 
\(
	\tilde{u}^\nu 
  {}={}
	(\tilde{\chi}^{\nu}, \tilde{\psi}^{\nu}, \tilde{\tau}^{\nu})
\).
We compute 
\(
	\bar{\varsigma}^\nu
{}={}
	\bar{\psi}^{\nu}
{}-{}
	2\tilde{\psi}^{\nu}
{}+{}
	\psi^\nu.
\)
Let us also define the triplet 
\(
	(
	  \bar{x}^\nu, \bar{y}^\nu, \bar{s}^\nu
	)
{}\dfn{}
	(
	  \bar{\chi}^{\nu}/\bar{\tau}^{\nu}, 
	  \bar{\psi}^{\nu}/\bar{\tau}^{\nu}, 
	  \bar{\varsigma}^{\nu}/\bar{\tau}^{\nu}
	)
\), which serves as the candidate primal-dual solution at iteration \(\nu\).
The relative primal residual is
\begin{equation}
 \mathrm{pr}_{\nu} {}={} \frac
		      {
			  \|
			  A\bar{x}^\nu + \bar{s}^\nu - b
			  \|
		      }
		      {
			  1+\|b\|
		      }
\end{equation}
The relative dual residual is 
\begin{equation}
 \mathrm{dr}_\nu {}={} \frac
		      {
			  \|
			  A^{\tttop}\bar{y}^\nu + c
			  \|
		      }
		      {
			  1+\|c\|
		      }
\end{equation}
The relative duality gap is defined as 
\begin{equation}
 \mathrm{gap}_\nu {}={} \frac
		      {
			  |\<c, \bar{x}^{\nu}\> {}+{} \<b,\bar{y}^{\nu}\>|
		      }
		      {
			  1 + |\<c, \bar{x}^{\nu}\>| + |\<b, \bar{y}^{\nu}\>|
		      }
\end{equation}
If \(\mathrm{pr}_{\nu}\), \(\mathrm{dr}_{\nu}\) and \(\mathrm{gap}_{\nu}\) are all below a specified tolerance
\(\epsilon>0\), then we conclude  that~\eqref{eq:cone-program-1} is feasible, the algorithm is terminated
and the triplet 
\(
	(
	  x^\nu, y^\nu, s^\nu
	)
\)
is an approximate solution.

The relative infeasibility certificate is defined as
(note that \(\bar{y}^{\nu}\in \K^*\) as a result of the projection step)
\begin{equation}
 \mathrm{ic}_{\nu} {}={}
	      \begin{cases}
		      {
			 \left\|b\right\| 
			 \left\|A^{\tttop}\bar{y}^{\nu}\right\|
		      }
		      /
		      {
			  \<b, \bar{y}^{\nu}\>
		      },
		      &\text{ if } \<b, \bar{y}^{\nu}\>{}<{}0\\
		      +\infty,&\text{ else}
	      \end{cases}
\end{equation}
Likewise, the relative unboundedness certificate is defined as
\begin{equation}
 \mathrm{uc}_{\nu} {}={}
	      \begin{cases}		     
		      {
			  \left\|c\right\| 
			  \left\|A\bar{x}^{\nu} + \bar{s}^{\nu}\right\|
		      }
		      /
		      {
			  \<c, \bar{x}^{\nu}\>
		      },
		      &\text{ if } \<c, \bar{x}^{\nu}\> {}<{} 0\\
		      +\infty,&\text{ else}
	      \end{cases}
\end{equation}
Provided that \(\bar{u}^{\nu}\) is not a feasible \(\epsilon\)-optimal point, it is a certificate of
unboundedness if \(\mathrm{uc}_{\nu}<\epsilon\) and it is a certificate of infeasibility if 
\(\mathrm{ic}_{\nu}<\epsilon\).
\end{subequations}

\subsection{Quasi-Newtonian directions}\label{sec:quasi-newtonian}
Quasi-Newtonian directions \(d^{\nu}\) can be computed according to the general rule
\begin{equation}
	d^\nu
{}={}
	-B_{\nu}^{-1} Ru^{\nu}
{}={}
	-H_{\nu} Ru^{\nu},
\end{equation}
where invertible linear operators \(H_{\nu}\) are updated according to certain 
low-rank updates so as to satisfy certain secant conditions starting 
from an initial operator \(H_0\).

\subsubsection{Restarted Broyden directions}
Here we make use of Powell's trick to update linear operators \(B_\nu\)
in such a way so as to enforce nonsignularity using the recursive 
formula~\cite{powellBroyden,supermann}:
\begin{equation}
	B_{\nu+1}
{}={}
	B_{\nu}
{}+{}
	\tfrac{1}{\|z^\nu\|^2}
	(\tilde{\xi}^\nu - B_{\nu}z^\nu)z^{\nu\top}
\end{equation}
where 
\(
	z^{\nu}
{}={}
	w^{\nu} 
{}-{}
	u^{\nu}
\),
\(
	\xi^{\nu}
{}={}
	Rw^{\nu} 
{}-{}
	Ru^{\nu}
\)
and for a fixed parameter \(\bar{\vartheta} \in (0,1)\) and 
\(
	\gamma_\nu 
{}={} 
	\<H_\nu \xi^{\nu},  z^{\nu}\>/\|z^\nu\|^2
\)
we have 
\(
	\tilde{\xi}^{\nu}
{}={}
	(1-\theta_\nu)B_{\nu}z^\nu + \theta_{\nu} \xi^{\nu}
\)
with
\begin{equation}\label{eq:theta-rbroyden}
	\theta_{\nu}
{}={}
	\begin{cases}
	 1		& \text{if } |\gamma_{\nu}| \geq \bar{\vartheta}
	 \\
	 \frac
	  {1-\operatorname{sgn}(\gamma_\nu)\bar{\vartheta}}
	  {1-\gamma_\nu}
		        & \text{otherwise}
	\end{cases}
\end{equation}
and the convention \(\operatorname{sgn}(0)=1\).
Using the Sherman-Morisson formula, operators \(H_{\nu}\) are updated as
follows:
\begin{equation}
	H_{\nu+1}
{}={}
	H_{\nu}
{}+{}
	\tfrac{1}{\< H_{\nu}\tilde{\xi}^{\nu}, z^{\nu}\>}  
			      (
			      z^{\nu} 
			{}-{} 
			      H_{\nu}\tilde{\xi}^{\nu}
			      )
			      (
				z^{\nu\top} H_{\nu}
			      ),
\end{equation}
thus lifting the need to compute and store matrices \(B_\nu\).

The Broyden method requires that we store matrices of dimension \((m{}+n{}+1) {}\times{} (m{}+n{}+1)\).
Here we employ a limited-memory \ac{rb} method which affords us a computationally
favorable implementation using buffers of length \(\mathrm{mem}\), that is
\(
	Z_{\nu}
{}={}
	[
	  z^\nu 
	  ~ 
	  z^{\nu-1} 
	  ~ 
	  \cdots 
	  ~ 
	  z^{\nu-\mathrm{mem}+1} 
	],
\) and 
\(
	\tilde{Z}_{\nu}
{}={}
	[
	  \tilde{z}^\nu 
	  ~ 
	  \tilde{z}^{\nu-1} 
	  ~ 
	  \cdots 
	  ~ 
	  \tilde{z}^{\nu-\mathrm{mem}+1} 
	],
\)
where \(\tilde{z}^{i}\) are the auxiliary variables 
\(
	\tilde{z}^{i} 
{}\dfn{}
	\nicefrac{z^i {}-{} H_i\tilde{\xi}^i}{\<s^i, H_i\tilde{\xi}^i\>}
\).
\begin{algorithm} 
\footnotesize
\begin{algorithmic}
 \REQUIRE Old buffers \(Z{}={}Z_{\nu}\) and \(\tilde{Z}{}={}\tilde{Z}_{\nu}\), 
	  \(\xi{}={}\xi^{\nu}\), 
	  \(r{}={}Ru^{\nu}\),
	  \(z{}={}z^\nu\),
	  \(\bar{\vartheta}\),
	  \(\mathrm{mem}\)
 \ENSURE  Direction \(d\), 
	  New buffers
 \STATE \(d {}\gets{} -r\), \(\tilde{z} {}\gets{} \xi\), \(m' {}\gets{}\) current cursor position
 \FOR {\(i=1,\ldots,m'\)}
    \STATE 
	\(
		\tilde{z} 
	{}\gets{}
		\tilde{z}
	{}+{}
	       \<z^{i},\tilde{z}\>\tilde{z}^i
	\), and
	\(
		d
	{}\gets{}
		d
	{}+{}
	       \<z^{i},d\>\tilde{z}^i
	\)
 \ENDFOR
 \STATE Compute \(\theta\) as in~\eqref{eq:theta-rbroyden} with 
	\(
		\gamma 
	{}={} 
		\nicefrac{\<\tilde{z}, z\>}{\|z\|^2}
	\)
 \STATE \(
		\tilde{z}
	{}\gets{}
		(1-\theta)z
	{}+{}
		\theta\tilde{z}
	\), 
	\(
		\tilde{z}
	{}\gets{}
		\nicefrac{z-\tilde{z}}{\<z,\tilde{z}\>}
	\), and
	\(
		d
	{}\gets{}
		d
	{}+{}
		\<z,d\>\tilde{z}
	\)
 \IF {\(m'=\mathrm{mem}\)}
    \STATE Empty buffers \(Z\) and \(\tilde{Z}\), \(m'\gets 1\)
 \ELSE
    \STATE Append \(z\) to \(Z\) and \(\tilde{z}\) to \(\tilde{Z}\), \(m'\gets m'+1\)
 \ENDIF
\end{algorithmic}
\caption{Modified restarted Broyden method}\label{alg:restarted-broyden}
\end{algorithm}
We have observed that SuperSCS performs better when deactivating K0 steps or using a 
small value for \(c_0\) (e.g., \(c_0=0.1\)) when using restarted Broyden directions.

\subsubsection{Anderson's acceleration}
\ac{aa} imposes a multi-secant condition~\cite{WalkerNiAnderson2011,Fang2009}.
In particular, at every iteration \(\nu\) we update a buffer of \(\mathrm{mem}\)
past values of \(z\) and \(\xi\), that is we construct a buffer \(Z_{\nu}\) as 
above and a buffer
\(
	\Xi_\nu 
{}={}
	\begin{bmatrix}
	  \xi^\nu 
	  & 
	  \xi^{\nu-1} 
	  & 
	  \cdots 
	  & 
	  \xi^{\nu-{\mathrm{mem}}+1} 
	\end{bmatrix}
\).
Directions are computed as 
\begin{equation}
	  d^\nu 
{}={}
	  -Ru^\nu 
{}-{} 
	  (Z_{\nu}{}-{}\Xi_{\nu})t^{\nu},
\end{equation}
where \(t^\nu\) is a least-squares solution of the linear system
\(\Xi_\nu t^\nu  = Ru^\nu,\) that is \(t^{\nu}\) solves
\begin{equation}
 \minimize_{t^\nu}\|\Xi_\nu t^\nu  - Ru^\nu\|^2,
\end{equation}
and can be solved using the singular value decomposition of \(\Xi_{\nu}\),
or a QR factorization which may be updated at every iteration~\cite{WalkerNiAnderson2011}.
In practice Anderson's acceleration works well for short memory lengths,
typically between 3 and 10, and, more often than not, outperforms 
the above restarted Broyden directions.

\subsection{Convergence}
\begin{figure}[h]
 \centering
 \includegraphics[width=0.99\linewidth]{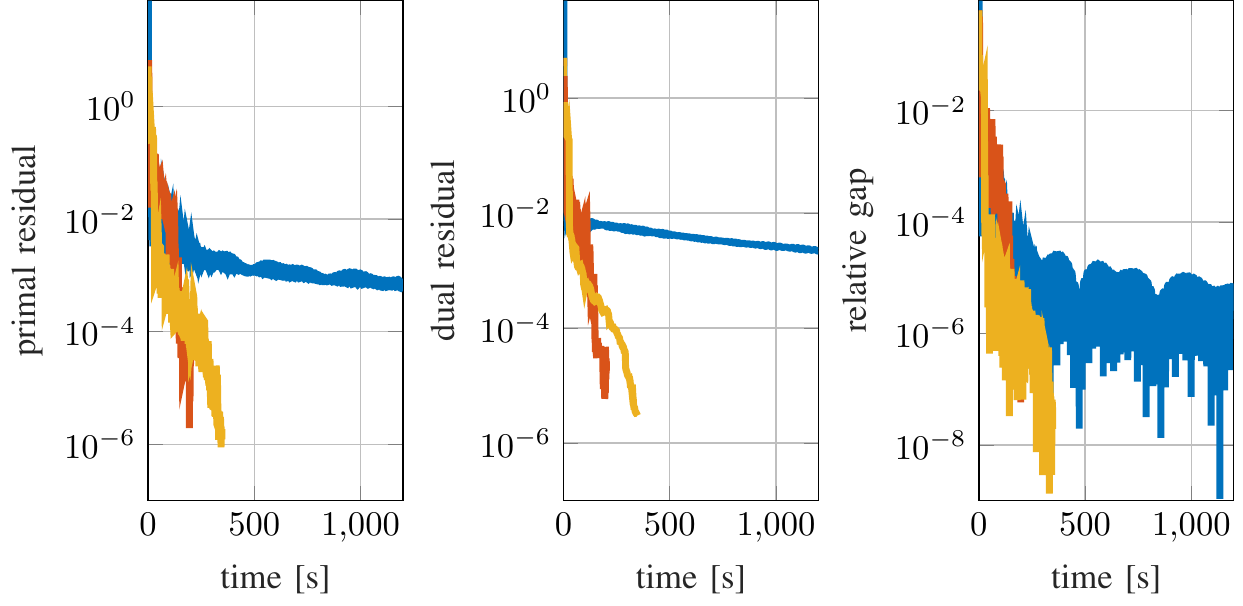} 
 \definecolor{mycolor1}{rgb}{0.00000,0.44700,0.74100}%
 \definecolor{mycolor2}{rgb}{0.85000,0.32500,0.09800}%
 \definecolor{mycolor3}{rgb}{0.92900,0.69400,0.12500}%
 \caption{Progress of SuperSCS (with RB and AA directions) and SCS versus time for a large-scale
          SDP of the form \eqref{eq:rpca} with \(d=500\) (with \(m=625751\) and \(n=250501\)).
          [\textcolor{mycolor1}{\bf ---} SCS; 
           \textcolor{mycolor2}{\bf ---} SuperSCS RB with memory 50;
           \textcolor{mycolor3}{\bf ---} SuperSCS AA with memory 5].}
 \label{fig:progress_superscs_vs_scs}           
\end{figure}

The convergence properties of SuperSCS are inherited by those of the 
general SuperMann scheme \cite{supermann}. In particular, under a weak 
 boundedness assumption for the quasi-Newtonian directions \(d^{\nu}\),
\(Ru^{\nu}\) converges to zero and \(u^{\nu}\) converges to a \(u^\star\)
which satisfies the \ac{hsde} \eqref{eq:monotone_inclusion}. If, additionally,
\(R\) is metrically subregular at \(u^{\star}\) --- a weak assumption ---
then, the convergence is R-linear. Under additional assumptions, SuperSCS
with the full-memory counterpart of the above restarted Broyden 
scheme can be shown to converge superlinearly. The restarted Broyden 
directions of Algorithm \ref{alg:restarted-broyden} and Anderson's 
acceleration, though not proven superlinear directions, exhibit steep 
linear convergence as shown in the next section (see Fig. \ref{fig:progress_superscs_vs_scs})
and have low memory requirements. 
\section{Benchmarks and results}\label{sec:benchmarks}
\subsection{Benchmarking methodology}
In order to compare different solvers in a statistically meaningful way, 
we use the \ac{dm} plot~\cite{dolanMore} and the shifted geometric mean%
~\cite{mittelmanSGM}. 
The \ac{dm} plot allows us to 
compare solvers in terms of their relative performance (e.g., computation time, flops, 
etc) and robustness, i.e., their ability to successful solve a given problem up to a
certain tolerance.

Let \(P\) be a finite set of test problems and \(S\) a finite set of solvers we want to 
compare to one another. Let \(t_{p,s}\) denote the computation time that solver \(s\)
needs to solve problem \(p\). We define the ratio between \(t_{p,s}\) and the lowest
observed cost to solve this problem using a solver from \(S\) as 
\begin{equation}\label{eq:performance_ratio}
	r_{p,s} 
{}={}
	\tfrac{t_{p,s}}{\min_{s'\in S}t_{p,s'}}.
\end{equation}
If \(s\) cannot solve \(p\) at all, we define \(t_{p,s}=+\infty\) and 
\(r_{p,s}=+\infty\). The cumulative distribution of the performance ratio is the 
\ac{dm} performance profile plot. In particular, define
\begin{equation}
	  \rho_s(\tau) 
{}={}
	  \nicefrac{1}{|P|} {}\cdot{} |\{p\in P: r_{p,s}\leq \tau\}|,
\end{equation}
for \(\tau{}\geq{}1\). 
The \ac{dm} performance profile plot is the plot of \(\rho_s(\tau)\) versus 
\(\tau\) on a logarithmic \(x\)-axis. For every \(s\in S\), the value \(\rho_s(1)\)
is the probability of solver \(s\) to solve a given problem faster than all 
other solvers, while \(\lim_{\tau\to\infty}\rho_s(\tau)\) is the probability that solver \(s\) solves 
a given problem at all.

As demonstrated in~\cite{gould2016}, \ac{dm} plots aim at comparing multiple 
solvers to the best one (cf.~\eqref{eq:performance_ratio}), than to one another.
Therefore, alongside we shall report
the shifted geometric mean of computation times for each solver following~\cite{mittelmanSGM}. 
For solver \(s{}\in{}S\), define the vector \(t^s{}\in{}\R^{|P|}\) with
\[
 t^s_p = \begin{cases}
          t_{p,s}&\text{ if } t_{p,s}{}<{}\infty
          \\
          100\max\{t_{p,s}{}\mid{}p\in P, t_{p,s}{}<{}\infty\} & \text{ otherwise}
         \end{cases}
\]
The shifted geometric mean of \(t^{s}\) with shifting parameter \(\sigma{}\geq{}0\) 
is defined as
\begin{equation}
	\operatorname{sgm}_{\sigma}
{}={} 
	\exp
	\Big[
	  \textstyle\sum_{p{}\in{}P}\ln
	      \left(
		  \max\{1, \sigma {}+{} t^{s}_{p}\}
	      \right)
	\Big]
{}-{}
	\sigma.
\end{equation}
Hereafter we use \(\sigma = \unit[10]{s}\).

In what follows we compare SuperSCS with the quasi-Newtonian direction methods presented
in Section~\ref{sec:quasi-newtonian} against SCS~\cite{scs,ODonoghue2016}. All tolerances
are fixed to \(10^{-4}\). In order to allow for a fair comparison among algorithms with 
different per-iteration cost, we do not impose a maximum number of iterations; instead, 
we consider that an algorithm has failed to produce a solution --- or a certificate
of unboundedness/infeasibility --- if it has not terminated after a certain (large) maximum time.
All benchmarks were executed on a system with a quad-core i5-6200U CPU at \(\unit[2.30]{GHz}\) and 
\(\unit[12]{GB}\) RAM running Ubuntu 14.04.

\subsection{Semidefinite programming problems}
\begin{figure}[H]
    \centering
    \includegraphics[width=0.8\columnwidth]{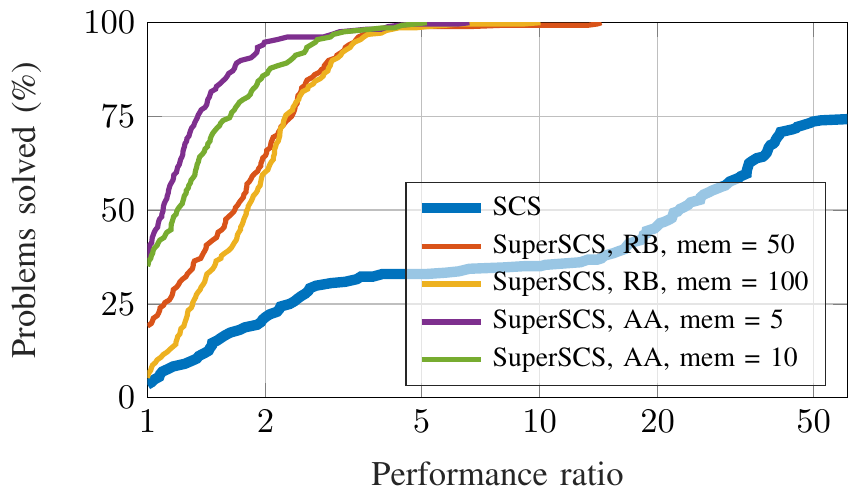} 
    \caption{\ac{dm} performance plot on 288 \(\ell_1\)-regularized PCA problems of the form~\eqref{eq:rpca}.}
    \label{fig:reg-pca-dmp}
\end{figure}

Let us consider the problem of 
sparse principal component analysis with an \(\ell_1\)-regularizer which has the form~\cite{sparsePCA}.
\begin{subequations}
\label{eq:rpca}
\begin{align}
 \maximize&\	   \operatorname{trace}(SZ) - \lambda\|Z\|_1
\\
 \stt&\ 
	  \operatorname{trace}(Z) = 1,\ Z=Z^\top,\ Z\succeq 0
\end{align}
\end{subequations}
A total of 288 randomly generated problems was used
\begin{wraptable}{l}{0.27\textwidth}
\scriptsize
\captionof{table}{Regularized PCA SDP: Solver Statistics}
\vspace{-5pt}
\label{tab:reg-pca}
\begin{center}
  \begin{tabular}{l|l|l}
    \hline
    \textbf{Method}         &   \(\operatorname{sgm}_{10}(t^{s})\) & Success
    \\ \hline 
    SCS                     &   96.82                        & 74.31\%
    \\ 
    \ac{rb} (mem:50)  &   3.17                         & 100\%
    \\ 
    \ac{rb} (mem:100) &   3.31                         & 100\%
    \\
    \ac{aa} (mem:5)   &   2.13                         & 100\%
    \\
    \ac{aa} (mem:10)  &   2.66                         & 100\%
    \\
    \hline
  \end{tabular}
\end{center}
\vspace{-10pt}
\end{wraptable} for benchmarking with 
\(d \in \{50, 120,$ $140, 180\}\) and \(\lambda \in $ $\{0.1,2, 5\}\), where $d$ is the dimension of $Z$. 
The DP plot in Fig. \ref{fig:reg-pca-dmp} shows that SuperSCS is consistently
faster and more robust compared to SCS. In Table \ref{tab:reg-pca} we see that
SuperSCS is faster than SCS by more than an order of magnitude; in particular,
SuperSCS with AA directions and memory \(5\) was found to perform best.
%
%

Additionally, for this benchmark we tried the interior-point solvers of SDPT3~\cite{sdpt3} 
and Sedumi~\cite{S98guide}. Both exhibited similar performance; in particular, for problems
of dimension $Z\in\R^{140\times 140}$, SDPT3 and Sedumi required \(\unit[6500]{s}\) to \(\unit[7500]{s}\)
and for problems of dimension \(Z\in\R^{180\times 180}\) they required \(\unit[19000]{s}\) 
to \(\unit[21500]{s}\). Note that SuperSCS (with RB and memory 50) solves all 
problems in no more than \(\unit[11.7]{s}\). Additionally, interior-point methods 
have an immense memory footprint of several GB, in this example whereas SuperSCS with RB 
and memory 100 needs as little as \(\unit[211.1]{MB}\) and with AA and memory 5
consumes just \(\unit[46.2]{MB}\).

\subsection{LASSO problems}\label{sec:lasso-dmp}
Regularized least-squares problems with the \(\|{}\cdot{}\|_1\)-regularizer, also known as 
LASSO problems, are optimization problems of the form
\begin{equation}
 \label{eq:lasso}
 \minimize_{x\in\R^n} \nicefrac{1}{2}\|Ax-b\|^2 + \mu \|x\|_1,
\end{equation}
where \(A\in\R^{m\times n}\) is a (sparse) matrix, and \(\mu>0\) is the regularization weight.
LASSO problems find applications in statistics and compressed 
sensing~\cite{recursiveCompressedSensing}. LASSO problems are cast as second-order cone 
programs~\cite{LOBO1998193}.

\begin{figure}[h]
      \centering
      \includegraphics[width=0.8\columnwidth]{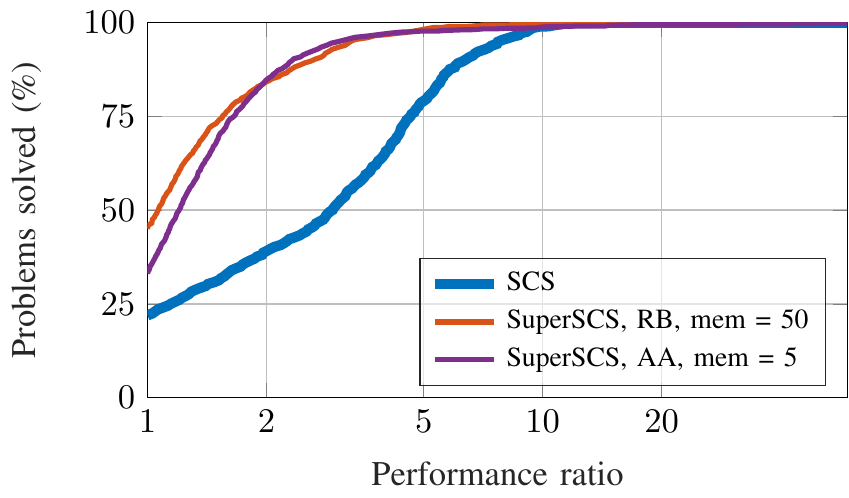} 
      \caption{\ac{dm} performance plot on 1152 LASSO problems.}
      \label{fig:lasso-dmp}
\end{figure}

\begin{wraptable}{l}{0.27\textwidth}
\vspace{-5pt}
\scriptsize
\captionof{table}{LASSO: Solver Statistics}
\vspace{-5pt}
\label{tab:lasso-sgm}
\begin{center}
  \begin{tabular}{ l | l | l}
    \hline
    \textbf{Method}         &   \(\operatorname{sgm}_{10}(t^{s})\) & Success
    \\ \hline 
    SCS                     &   5.97                               & 100\%
    \\ 
    \ac{rb} (mem:50)  &   2.88                               & 100\%
    \\ 
    \ac{rb} (mem:100) &   2.61                               & 100\%
    \\
    \ac{aa} (mem:5)   &   3.38                               & 100\%
    \\
    \ac{aa} (mem:10)  &   3.87                               & 100\%
    \\\hline
  \end{tabular}
\end{center}
\vspace{-15pt}
\end{wraptable}
LASSO problems aim at finding a sparse vector \(x\) which minimizes \(\|Ax-b\|^2\). The sparseness
of the minimizer \(x^\star\) can be controlled by \(\mu\in\{0.01, 0.1, 1\}\). 
We tested 1152 randomly generated LASSO problems with 
\(n\in\{631, 1000, 1585, 2512\}\),
\(m = \lceil n/5 \rceil\), and matrices \(A\) with condition numbers 
\(\kappa_A \in \{10, 215, 4600, 10^5\}\). 
The DM plot in Fig.~\ref{fig:lasso-dmp} and the statistics presented in 
Table \ref{tab:lasso-sgm} demonstrate that SuperSCS, both with RB and AA 
directions, outperforms SCS.

\subsection{Sparse \(\ell_1\)-regularized logistic regression}
\begin{figure}[H]
      \centering
      \includegraphics[width=0.8\columnwidth]{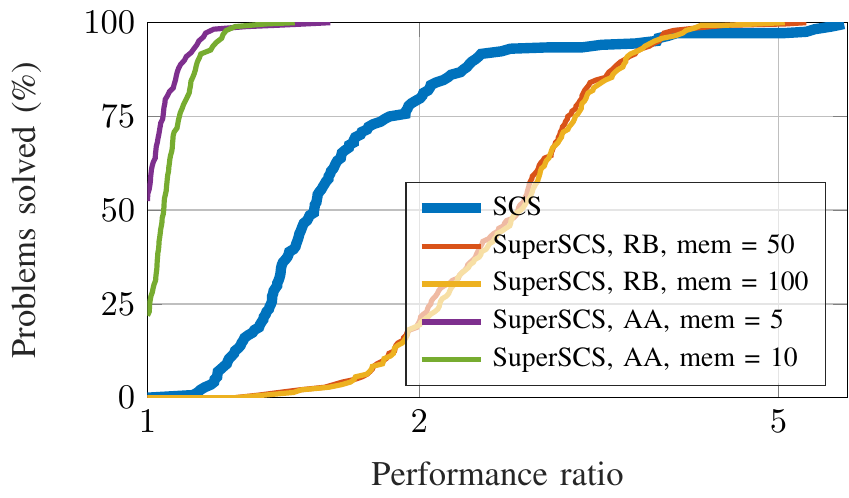} 
      \caption{\ac{dm} performance plot on 288 logistic regression problems of the form~\eqref{eq:logreg}.}
      \label{fig:logreg-dmp}
\end{figure}

Logistic regression is a regression model where dependent variables are binary~\cite{logregCox}.
The \(\ell_1\) regularized variant of logistic regression aims at performing simultaneous regression
and feature selection and amounts to solving an optimization problem of the following form~\cite{sparseLogReg}:
\begin{equation}
\label{eq:logreg}
	\minimize_{w\in\R^{p}} 
	\lambda \|w\|_1 
{}-{}
	\sum_{i=1}^{q}
	  \log(1 {}+{} \exp(a^\top w_i {}+{} b)).
\end{equation}
Similar to LASSO (Sec.~\ref{sec:lasso-dmp}), parameter \(\lambda>0\)
controls the sparseness of the solution \(w^\star\). Such problems can be cast as 
conic programs with the exponential cone, which is not self-dual.
A total of 288 randomly generated problems was used for benchmarking 
\begin{wraptable}{l}{0.265\textwidth}
\scriptsize
\captionof{table}{Sparse \(\ell_1\)-regularized logistic regression: Solver Statistics}
\vspace{-5pt}
\label{tab:logreg-dmp}
\begin{center}
  \begin{tabular}{ l | l | l}
    \hline
    \textbf{Method}         &   \(\operatorname{sgm}_{10}(t^{s})\) & Success
    \\ \hline 
    SCS                     &   4.85                         & 100\%
    \\ 
    \ac{rb} (mem:50)  &   7.23                         & 100\%
    \\ 
    \ac{rb} (mem:100) &   7.28                         & 100\%
    \\
    \ac{aa} (mem:5)   &   3.00                         & 100\%
    \\
    \ac{aa} (mem:10)  &   3.09                         & 100\%
    \\
    \hline
  \end{tabular}
\end{center}
\vspace{-12pt}
\end{wraptable} with 
\(p\in\{80, 100\}\), \(q\in\{50, 100, 120\}\) and \(\lambda \in \{10, 20, 50\}\).
In  Fig.~\ref{fig:logreg-dmp} we observe that SuperSCS with \ac{rb} 
directions is slower compared to SCS, however SuperSCS with \ac{aa} is noticeably faster.

%
%

\subsection{Maros-M{\'e}sz{\'a}ros QP problems}
\begin{figure}[h]
    \centering
    \includegraphics[width=0.8\columnwidth]{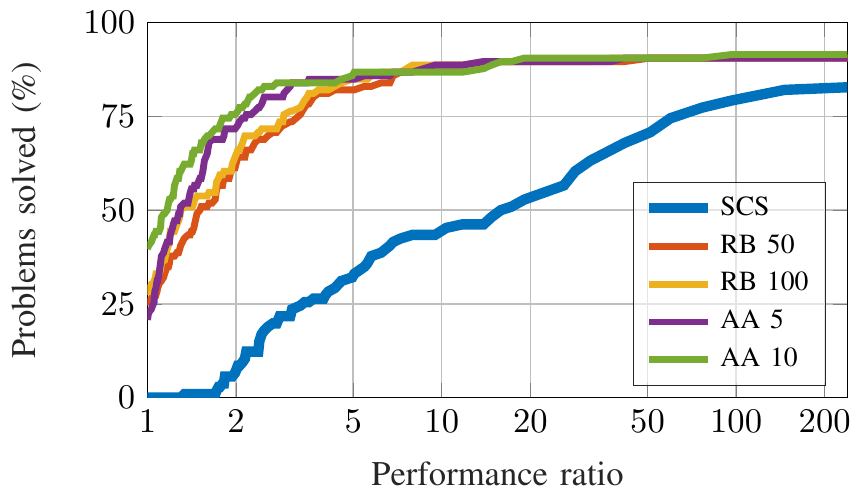} 
    \caption{\ac{dm} performance plot on the problems of the Maros-M{\'e}sz{\'a}ros repository.}
    \label{fig:maros-meszaros-dmp}
\end{figure}

\begin{wraptable}{l}{0.27\textwidth}
\vspace{-7pt}
\scriptsize
\captionof{table}{Maros-M{\'e}sz{\'a}ros QP problems: Solver Statistics}
\vspace{-9pt}
\begin{center}
  \begin{tabular}{ l | l | l}
    \hline
    \textbf{Method}      &  \(\operatorname{sgm}_{10}(t^{s})\)  & Success
    \\ \hline 
    SCS                           &   56.61                              & 83.02\%
    \\ 
    \ac{rb} (mem:50)  &   9.66                               & 90.57\%
    \\ 
    \ac{rb} (mem:100)  &   6.57                               & 90.57\%
    \\
    \ac{aa} (mem:5)    &   5.79                               & 90.57\%
    \\
    \ac{aa} (mem:10) &   8.62                               & 91.51\%  
    \\\hline
  \end{tabular}
\end{center}
\vspace{-13pt}
\end{wraptable}
Here we present performance of SCS and SuperSCS (with \ac{rb} and \ac{aa} directions)
on this collection of problems on the Maros-M{\'e}sz{\'a}ros collection of 
problems~\cite{Maros1999}.
As shown in Fig.~\ref{fig:maros-meszaros-dmp} SuperSCS, both with \ac{rb} directions
and Anderson's acceleration, is faster and more robust compared to SCS.
The two quasi-Newtonian directions appear to be on a par.
%
%

\section{Conclusions}
In this work we introduced SuperSCS: a first-order method for large-scale
conic optimization problems which combines the low iteration cost of SCS and the 
fast convergence of SuperMann.
We have compared SuperSCS with SCS on a broad collection of conic
optimization problems of practical interest. Using Dolan-Mor{\'e} plots and 
runtime statistics, we demonstrated that SuperSCS with Anderson's acceleration
is faster and more robust than SCS.

The C implementation of SuperSCS builds up on SCS and is a free open-source
software. SuperSCS can be interfaced from MATLAB, Python, can be invoked 
via CVX, CVXPy and YALMIP and is also available as a Docker image (see 
\url{https://kul-forbes.github.io/scs/}).

\bibliographystyle{plain}
\vspace{-1pt}
\bibliography{superScsBibliography}

\end{document}